\documentclass{amsart}
\usepackage[top=1.3in, left=1.3in, right=1.3in, bottom=1.3in]{geometry}
\usepackage{amsxtra,amssymb,multicol,graphicx}

\title{WHAT IS ...?}
\author{Author names}
\begin{document} 
\rightline{\Huge \bf WHAT IS ... A Blender?}
\vskip1cm
\rightline{\it \huge Ch. Bonatti, S. Crovisier, L. J. D\'\i az, A. Wilkinson}
\vskip1in


\noindent

A blender is a compact invariant set on which a  diffeomorphism has a certain behavior, which forces topologically ``thin" sets to intersect in a robust way, producing rich dynamics.  The term ``blender" describes its function: to blend together  stable and unstable manifolds.  Blenders have been used to construct diffeomorphisms with surprising properties and have played an important role in the classification of smooth dynamical systems.

One of the original applications of blenders is also one of the more striking.  A diffeomorphism $g$ of a compact manifold is {\em robustly transitive}
if there exists a point $x$ whose orbit $\{g^{n}(x): n\geq 0\}$ is dense in the manifold, and moreover this property persists when $g$ is slightly perturbed. Until the 1990's  there were no known robustly transitive diffeomorphisms in the isotopy class of the identity map on any manifold.
Bonatti and D\'iaz  (Ann. of Math., 1996\footnote{This is also where the term ``blender" was coined.}) used blenders to construct robustly transitive diffeomorphisms  as perturbations of the identity map on certain $3$-manifolds.

\begin{figure}[ht]
\begin{center}
\includegraphics[scale=0.28]{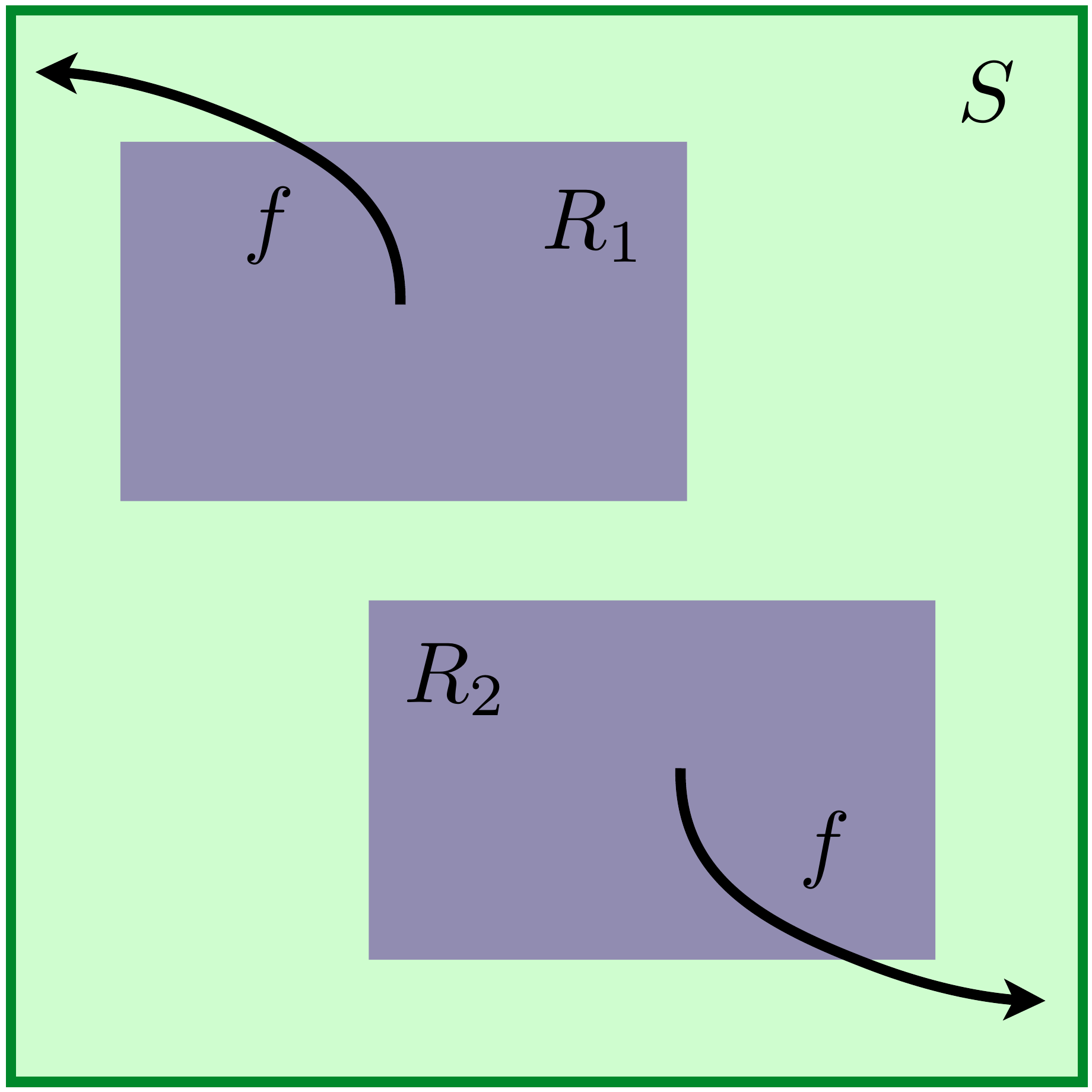}
\end{center}
\caption{An example of a proto-blender.   The map $f$ is defined on the union of the two rectangles $R_1$ and $R_2$ in the square $S$;
$f$ sends each $R_i$ onto the entire square  $S$ affinely, respecting the horizontal and vertical directions,  with the 
horizontal expansion factor  less than $2$. Note that $f$ fixes a unique point in each rectangle $R_i$.
 }
\end{figure}

To construct a blender one typically starts with a {\em proto-blender}; an example is the map 
$f$ pictured in Figure  1. 
The function $f$ maps each of the two rectangles $R_1$ and $R_2$ affinely onto the square $S$ and has the property that the vertical projections of $R_1$ and $R_2$ onto the horizontal direction overlap. Each rectangle contains a unique fixed point for $f$.

The compact set $\Omega = \bigcap_{n\geq 0} f^{-n} S$ is $f$-invariant, meaning $f(\Omega) = \Omega$, and
is characterized as the set of points in $S$ on which $f$ can be iterated infinitely many times:   $x\in \Omega$ if and only if  $f^n(x)\in S$ for all $n\geq 0$.
 $\Omega$ is a Cantor set, obtained by intersecting all preimages $f^{-i}(S)$ of the square, 
which nest in a regular pattern as in Figure 2.

\begin{figure}[ht]
\begin{center}\includegraphics[scale=0.35]{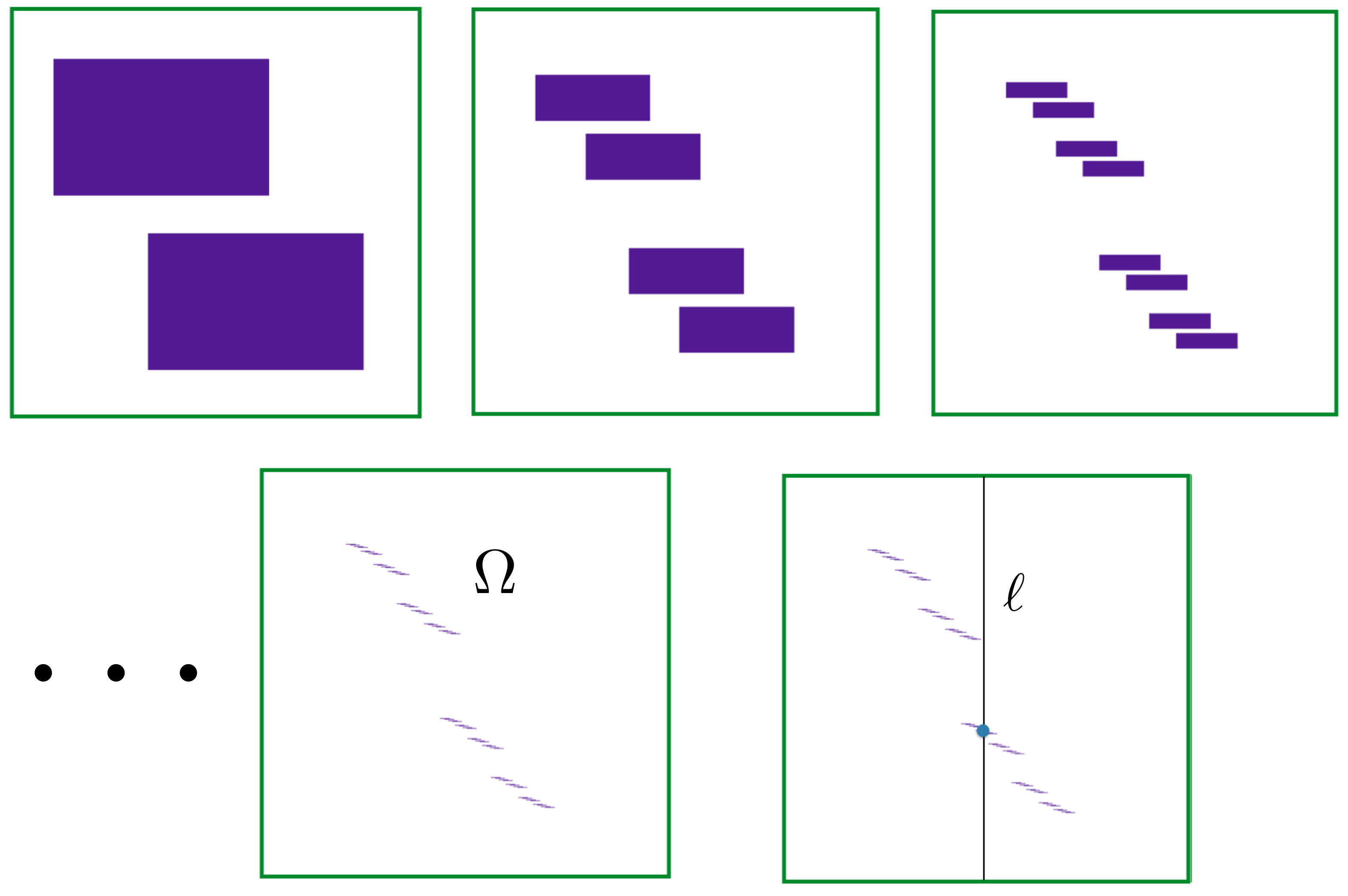}
\end{center}
\caption{The invariant Cantor set  $\Omega$ produced by the proto-blender $f$ is the nested intersection of preimages of $S$ under $f$.  Any vertical line segment $\ell$ close to the center of the square intersects $\Omega$ in at least one point.  The line segment can be replaced by segment with nearly vertical slope, or even a smooth curve nearly tangent to the vertical direction.}
\end{figure}

Any vertical line $\ell$ between the fixed points in $R_1$ and in $R_2$ will meet $\Omega$.  To prove this, it is enough to see that  for every $i$ the
vertical projection of the set $f^{-i}(S)$  (consisting of $2^i$ horizontal rectangles) onto the horizontal is an interval. 
This can be checked inductively, observing that the projection of $f^{-i-1}(S)$ is the union of two re-scaled copies of the projection of $f^{-i}(S)$, 
which overlap. 

A more careful inspection of this proof reveals that the intersection is {\em robust} in two senses: first,
the line $\ell$ can be replaced by a line whose slope is close to vertical, or even by a $C^1$ curve whose tangent vectors are close to vertical; secondly, the map $f$ can be replaced by any $C^1$ map $\hat f$ whose derivative is close to that of $f$.  Such an  $\hat f$
is called a {\em perturbation} of $f$.

The (topological) dimension of the Cantor set $\Omega$ is 0,   the dimension of $\ell$ is 1,  
the dimension of the square is $2$, and  $0+1 < 2$.  From a topological point of view, one would not expect these sets to intersect each other.  But from a {\em metric} point of view, the fractal set $\Omega$, when viewed along nearly vertical directions, appears to be $1$-dimensional, allowing $\Omega$ to intersect a vertical line, robustly.  If the rectangles $R_1$ and $R_2$ had disjoint projections, the proto-blender property would be destroyed.

\begin{figure}[ht]
\begin{center}
\includegraphics[scale=0.32]{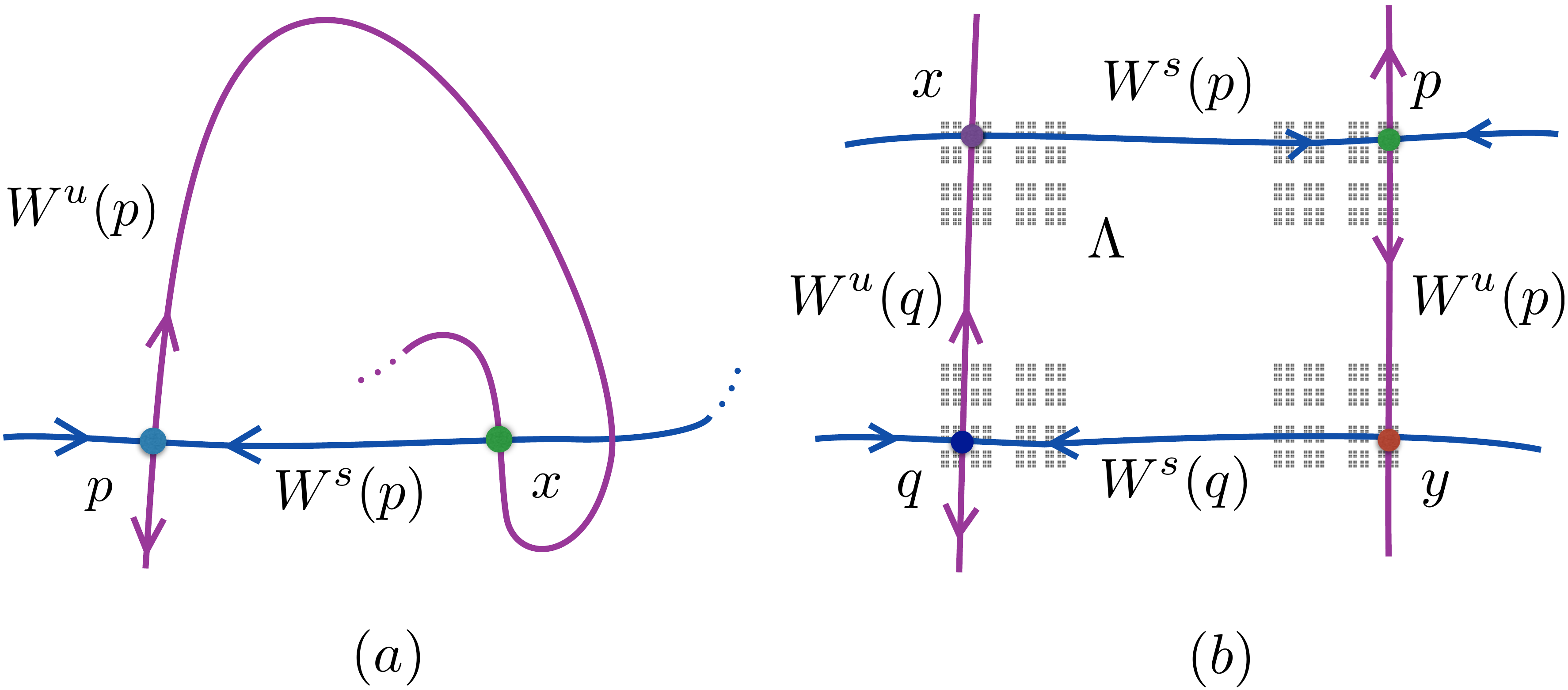}
\end{center}
\caption{(a) A transverse homoclinic intersection of stable and unstable manifolds, first discovered by Poincar\'e in his study of the three-body problem.  (b)  A horseshoe $\Lambda$ produced by a pair of transverse heteroclinic points $x$ and $y$.  Every point in the Cantor set $\Lambda$  can be approximated arbitrarily well both by a periodic point and by a point whose orbit is dense in $\Lambda$.   }
\end{figure}

This type of picture is embedded in a variety of smooth dynamical systems, where it is a robust mechanism for chaos.   The search for robust mechanisms for chaotic dynamics has a long history, tracing back to 
 Henri Poincar\'e's discovery of chaotic motion in the three-body problem of celestial mechanics.  Figure 3(a)  depicts the mechanism behind Poincar\'e's discovery, a local diffeomorphism of the plane with a saddle fixed point 
$p$ and another point $x$ whose orbit converges to $p$ both under forward and backward iterations (that is, under both the map and its inverse). 
Meeting at  $p$ are two smooth curves  $W^s(p)$ and  $W^u(p)$, the {\em stable and unstable manifolds} at $p$, respectively.
$W^s(p)$ is the set of points whose forward orbit  converges to $p$ and  $W^u(p)$ is the set of points whose backward   orbit converges to $p$.

In Figure 3(a),  the intersection of $W^s(p)$ and $W^u(p)$ is {\em transverse} at $x$: the tangent directions to 
$W^s(p)$ and $W^u(p)$ at $x$ span the set of all directions emanating from $x$ -- the tangent space at $x$ to the ambient manifold, in this case the plane.  The point $x$ is called a {\em transverse homoclinic point} for $p$. In Figure 3(b) a slight variation is depicted: here there are two periodic saddles $p$ and $q$ such that $W^s(p)$ and $W^u(q)$ intersect transversely at a point  $x$, and $W^u(p)$ and $W^s(q)$ intersect transversely at another point $y$.  The points $x$ and $y$ are called {\em transverse heteroclinic points} and they are arranged in a {\em transverse heteroclinic cycle}.

In the classification of the so-called {\em Axiom A diffeomorphisms}, carried out by Stephen Smale in the 1960's, 
transverse homoclinic and heteroclinic points play a central role.  Any transverse homoclinic point or heteroclinic cycle for a diffeomorphism is contained in in a special Cantor set $\Lambda$ called a {\em horseshoe}, an
invariant compact set  with strongly chaotic (or unpredictable) dynamical properties  (see [2] for a discussion).
Two notable properties of a horseshoe $\Lambda$ are:
\begin{enumerate}
\item Every point in $\Lambda$ can be approximated arbitrarily well by a periodic point in $\Lambda$.
\item There is a point in  $\Lambda$ whose orbit is dense in $\Lambda$.
\end{enumerate}

Horseshoes and periodic saddles are both examples of {\em hyperbolic sets}: a compact invariant set $\Lambda$ for a diffeomorphism $g$ is  hyperbolic if at every point in $\Lambda$  there are transverse stable and unstable manifolds $W^s(x)$ and $W^u(x)$ with $g(W^s(x)) = W^s(g(x))$ and $g(W^u(x)) = W^u(g(x))$.   For a large class
of diffeomorphisms known as Axiom A systems, Smale proved that the set of recurrent points can be decomposed into a disjoint union of  finitely many hyperbolic sets on which (1) and (2) above hold. 
This theory relies on the most basic property of transverse intersections, first investigated by Ren\'e Thom: robustness.  A transverse intersection of submanifolds cannot be destroyed by a small perturbation of the manifolds; in the dynamical setting, a transverse intersection of stable and unstable manifolds of two saddles cannot be destroyed by perturbing the diffeomorphism.

Classifying Axiom A systems was just the beginning.
To illustrate the limitations of the existing theory, Abraham and Smale constructed diffeomorphisms that are
{\em robustly non-hyperbolic}. These examples opened up the door for understanding a broader class of dynamics, and  blenders have turned out to be a key player in this emerging classification.

\begin{figure}[ht]
\begin{center}\includegraphics[scale=0.28]{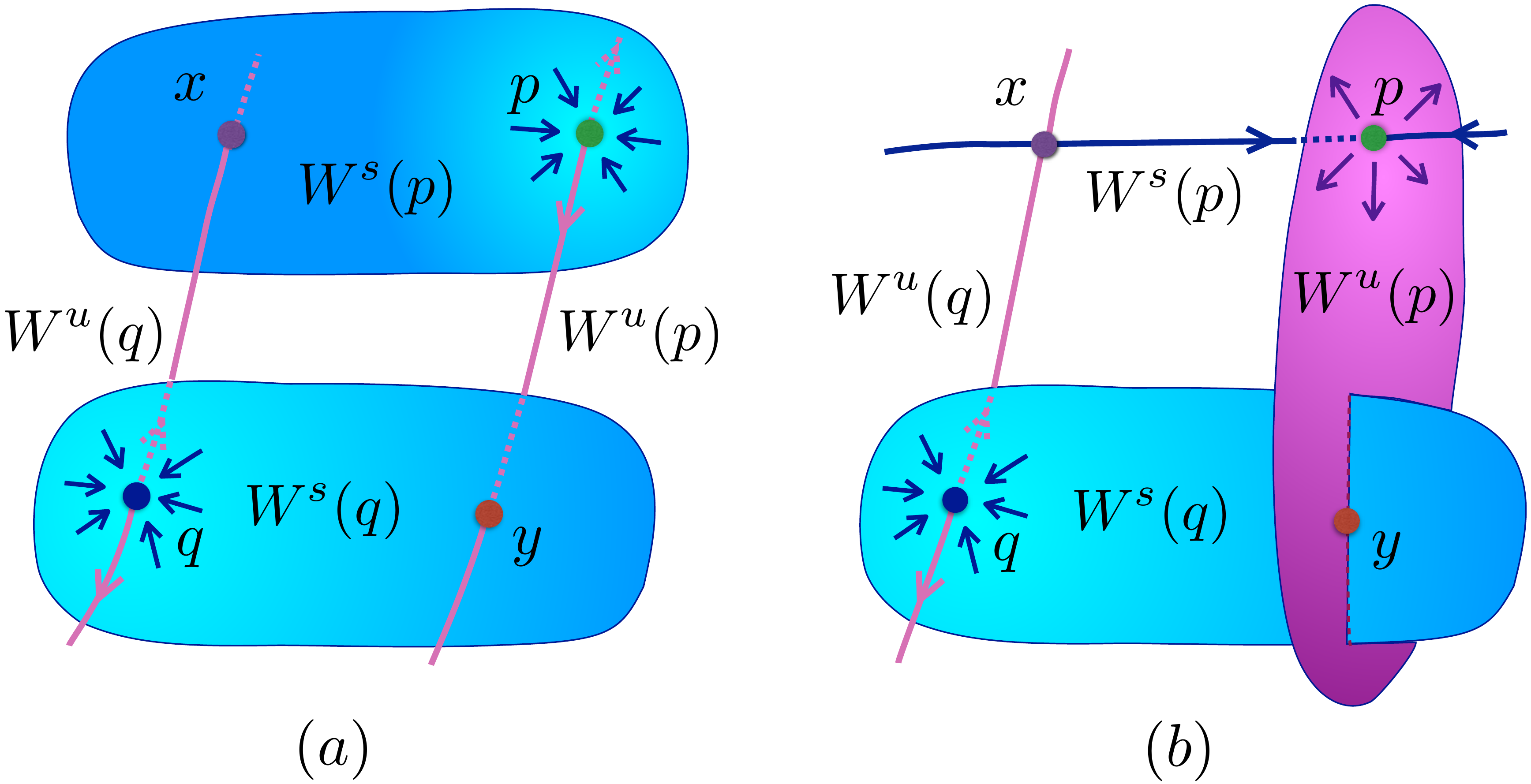}\end{center}
\caption{  (a) Transverse cycle.  (b) Non-transverse cycle.}
\end{figure}

Before constructing blenders and robust non-hyperbolic dynamics, we first illustrate (non-robust) dynamics of {\em non-hyperbolic} type.  To do so, let's return
to our example of two periodic saddles $p$ and $q$,  but this time in dimension $3$, where
saddle points can have stable and unstable manifolds of dimension either $1$ or $2$.
Suppose $p$ and $q$ are two  fixed points in dimension $3$ whose stable and unstable manifolds intersect.  If,
as in Figure 4(a),  the stable manifolds of 
$p$ and $q$ have the same dimension, then both intersections can be transverse\footnote{and the intersections are generically transverse in this case, a consequence of the Kupka-Smale Theorem.} producing  a horseshoe.

But quite another thing happens if the dimensions of the stable manifolds do not match up: the intersection between the $2$-dimensional
manifolds may be transverse,  but the other -- between $1$-dimensional manifolds  -- is {\em never} transverse, and thus cannot be robust. 
In the case depicted in figure 4(b), the orbit of the point $x$ accumulates on $q$ in the past and on $p$ in the future. The point $x$ cannot be contained in a hyperbolic set, because $W^s(p)$ and $W^s(q)$ have different dimensions.  On the other hand, this non-hyperbolicity is not robust, because this non-transverse intersection is easily destroyed by perturbation.

To obtain a robustly non-hyperbolic example,  we will replace the point $q$ in Figure 4(b) by a cube $Q$ containing a special type of horseshoe $\Lambda$ called a blender.   To produce $\Lambda$, we use the proto-blender $f\colon R_1\cup R_2\to S$  of Figure 1. 
The map $f$
has only expanding directions and is not injective; indeed, it  has precisely two inverse branches 
$f_1^{-1} \colon S\to R_1$ and  $f_2^{-1} \colon S\to R_2$.  In dimension three, we can embed 
these inverse branches into a local diffeomorphism by adding  a third,  expanded  direction, as detailed in Figure 5,
where the cube $Q$ is stretched and folded across itself by a local diffeomorphism $g$.

The horseshoe $\Lambda$ in Figure 5 is precisely the set of points whose orbits remain in the future and in the past in $Q$.
The set $W^u(\Lambda) $  of points in the cube that accumulate on $\Lambda$ in the past
is the cartesian product of the Cantor set $\Omega$ with segments parallel to the third, expanded direction.  $W^u(\Lambda)$ is the analogue of the unstable manifold of a saddle, but it is a fractal object rather than a smooth submanifold.

The set $\Lambda$ is an example of a blender,  and its main geometric property is that
 any vertical curve crossing $Q$ close enough the center intersects $W^u(\Lambda)$. 
 In other words, this blender is a horseshoe whose unstable set behaves likes a surface even though its topological dimension is one.  This property is robust.  While the definition of blender is still evolving as new constructions arise, a working definition is: {\em A blender is a compact hyperbolic set whose unstable set has dimension strictly less than one would predict by looking at its intersection with families of submanifolds.}

Figure 6 illustrates  robust nonhyperbolic dynamics, produced by combining Figure 4(b) with a blender.
The connection between the stable manifold of $p$ and  $W^{u}(\Lambda)$ cannot be destroyed by perturbation, and the  transverse  intersection between the unstable manifold of $p$ and the stable manifold of a point $z \in \Lambda$ is also robust.  The orbit of the point $z$ is contained in a compact invariant set with complicated dynamics, in particular satisfying property (1) above.

\begin{figure}[ht]
\begin{center}
\includegraphics[scale=0.28]{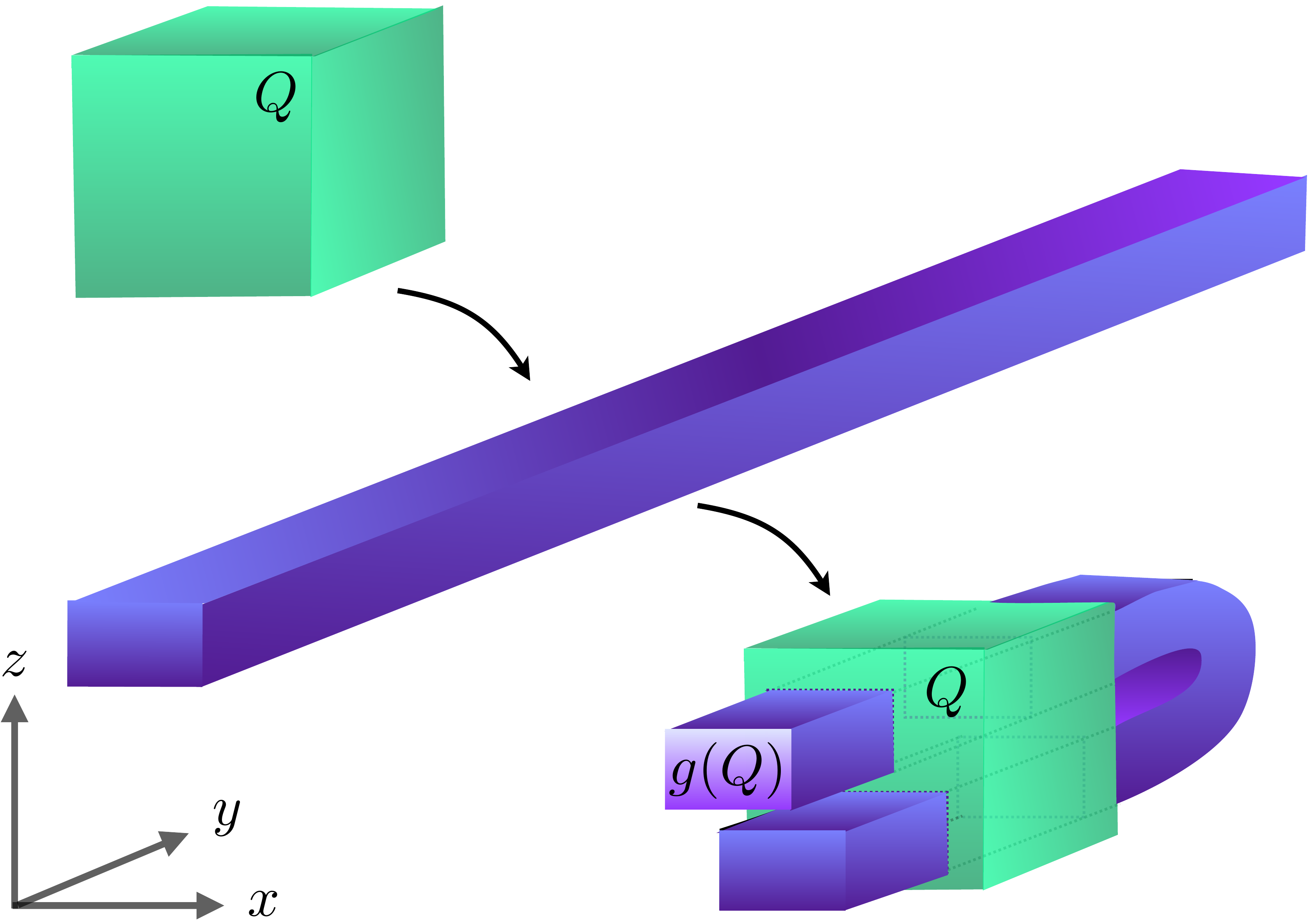}
\end{center}
\caption{Constructing a blender, a type of horseshoe with a proto-blender built into its contracting directions. In the cube $Q$ the local
 diffeomorphism $g$ contracts the segments in the axial directions parallel to the front face (the $xz$-plane), 
elongates the cube into the third axial direction (the $y$-axis), and then folds this elongated piece across the original cube $Q$, as pictured.   Each slice of
$Q \cap g(Q)$ parallel to the $xz$-plane resembles exactly the picture of $R_1\cup R_2$ in the square $S$.
The restriction of $g^{-1}$ to these rectangles in this slice just is a copy of the proto-blender $f$ from Figure 1,
whose image is another $xz$-slice of $Q$.}
\end{figure}

\begin{figure}[ht]
\begin{center}
\includegraphics[scale=0.28]{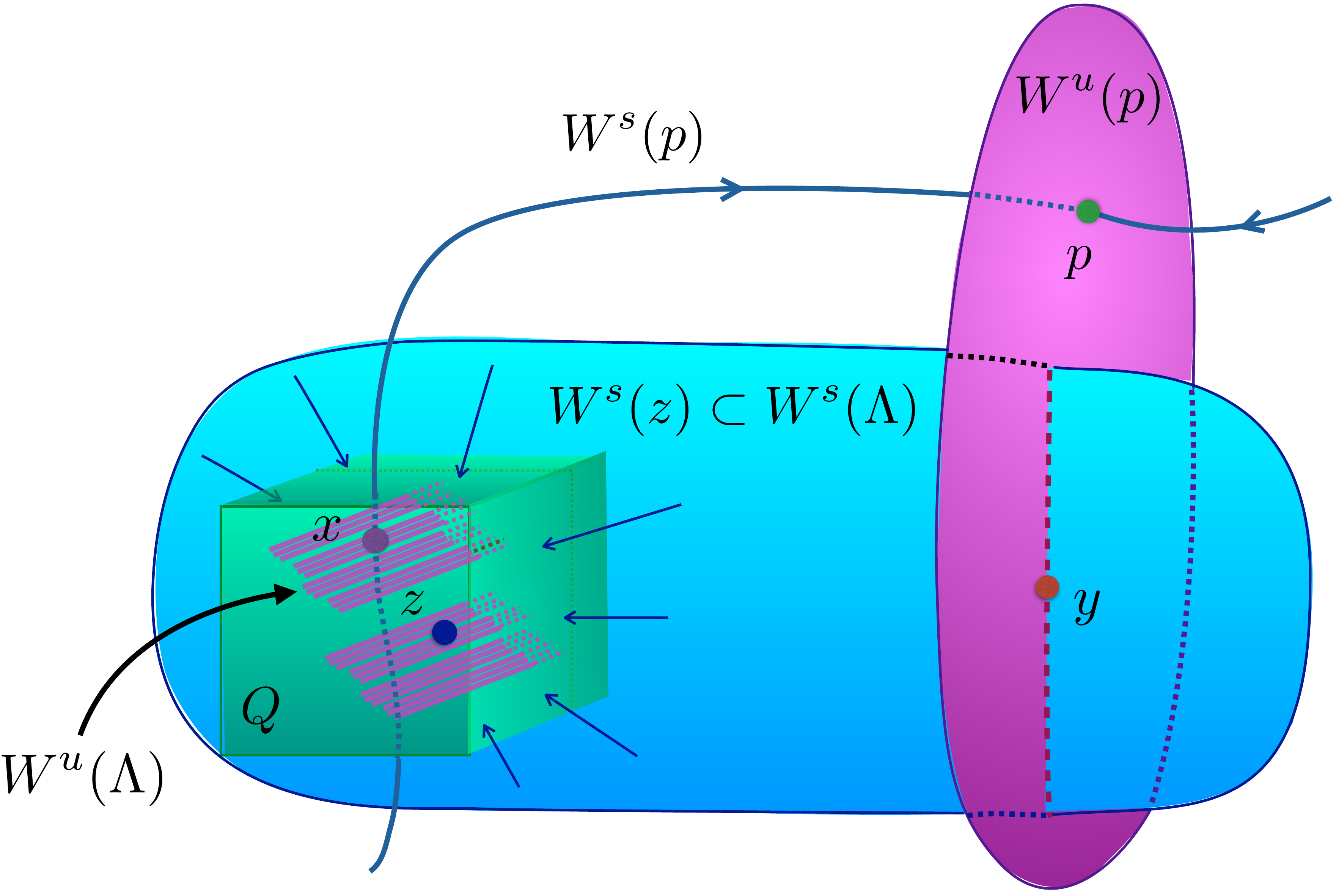}
\end{center}
\caption{Replacing the periodic point $q$ in Figure 4(b) with a cube $Q$ containing the blender of Figure 5.  The orbit of the points $x$ and $y$  accumulate both on the saddle $p$ and the blender horseshoe $\Lambda$, producing an invariant subset of the dynamics with complicated, non-hyperbolic dynamics.  } 
\end{figure}

 Blenders are not just a tool to produce robust 
non-hyperbolic dynamics, they are in fact one of the two conjectured mechanisms responsible for robust non-hyperbolicity, the other being homoclinic tangencies.
This is because, in contrast to the original Abraham-Smale construction, blenders appear in a natural 
way in local bifurcations.  Indeed, whenever a diffeomorphism has two saddles $p$ and $q$ 
with different stable dimensions and are dynamically related as in Figure 4(b), there is a perturbation 
that produces a blender.

\section*{\bf Further Reading}

\noindent

\noindent
[1] Christian Bonatti, Lorenzo D\'\i az and Marcelo Viana.
Dynamics beyond Uniform Hyperbolicity: A Global
Geometric and Probabilistic Approach,
Encyclopedia Math. Sci., Springer, 2004.

\noindent
[2] Michael Shub, ``What is... a horseshoe?" Notices Amer. Math. Soc. {\bf 52} (2005), 516--517.


\end{document}